\documentclass[12pt,a4paper,reqno]{amsart}
\usepackage{amssymb, amscd}

\newtheorem{thm}{Theorem}[section]

\newtheorem*{thm*}{Theorem}

\newtheorem{lem}[thm]{Lemma}
\newtheorem{prop}[thm]{Proposition}

\newtheorem{defn}[thm]{Definition}

\newtheorem*{remark*}{Remarks}

\newtheorem*{defn*}{Definition}
\newtheorem*{claim*}{Claim}

\theoremstyle{definition}
\newtheorem{remark}[thm]{Remark}

 \renewcommand{\sectionmark}[1]{}

\begin{document}

\title[Symmetric homogeneous diophantine equations of odd degree]
{Symmetric homogeneous diophantine equations of odd degree}
\author[M. A. Reynya]{M. A. Reynya}
\address{Ahva St. 15/14, Haifa, Israel}

\email{misha\_371@mail.ru}


\begin{abstract}
We find a parametric solution of an arbitrary symmetric
homogeneous diophantine equation of 5th degree in 6 variables
using two primitive solutions. We then generalize this approach to
symmetric forms of any odd degree by proving the following
results.

(1) Every symmetric form of odd degree $n\ge 5$ in $6 \cdot
2^{n-5}$ variables has a rational parametric solution depending on
$2n-8$ parameters.

(2) Let $F(x_1, \dots ,x_N)$ be a symmetric form of odd degree
$n\ge 5$ in $N=6 \cdot 2^{n-4}$ variables, and let $q$ be any
rational number. Then the equation $F(x_i)=q$ has a rational
parametric solution depending on $2n-6$ parameters.

The latter result can be viewed as a solution of a problem of
Waring type for this class of forms.
\end{abstract}

\maketitle

\baselineskip 20pt
\section {Introduction} \label{sec:intro}

Given any form $F(x_1, \dots ,x_N)$ with rational coefficients,
one can look for conditions under which the equation $F(x_1, \dots
,x_N)=0$ has a nonzero rational solution. It is known that if $N$
is big enough with respect to the degree of $F$ (which is assumed
odd), such a solution always exists (see \cite{Bi1}, \cite{Bi2};
more recent papers \cite{W1}, \cite{W2} provide explicit estimates
for $N$). In the present paper, we consider {\it symmetric} forms
$F(x_1, \dots ,x_N)$, i.e. forms invariant under the natural
action of the symmetric group $S_N$ on variables. For such a form
the existence of a nonzero rational solution is guaranteed, and we
are interested in getting {\it parametric} solutions and
estimating the number of parameters. More precisely, for any given
integers $n$ and $k$ we want to find an integer $N$ such that
every symmetric form in $N$ variables of degree $2n+1$ with
rational coefficients has a parametric solution depending on $k$
parameters.

Our method, being essentially elementary, is inspired by a
geometric approach proposed by Manin \cite{M} for studying
rational points on cubic hypersurfaces: given a few rational
points on such a variety, one can use a composition law in order
to produce new rational points. We implement a similar idea in the
case of symmetric forms of odd degree: starting from two
``primitive'' rational points, we propose an algorithm producing
lots of rational points.

To explain our method, we begin with the case of quintics in 6
variables (see Section \ref{sec:quin}). We prove that any such
quintic has a parametric solution $(x_1,\dots ,x_6)$ such that
$x_1+\dots +x_6=0$ (Theorem \ref{th:quintic}). Note that to the
best of our knowledge, such a result did not appear in the
literature, even in the simplest case of diagonal quintics. Our
method is significantly different both from the elementary
approach in \cite{R}, \cite{C} and computer investigation of
\cite{LP}, \cite{LPS}, \cite{E}. It can be compared with
geometric constructions in \cite{Br} which led to some new insight
on parametric solutions obtained in \cite{SwD}.

Next, in Section \ref{sec:gen}, we extend our method to symmetric
forms of arbitrary odd degree. Our main result (Theorem
\ref{th:general}) states that every symmetric form of odd degree
$n\ge 5$ in $6 \cdot 2^{n-5}$ variables has a rational parametric
solution depending on $2n-8$ parameters. As an application, we
obtain a solution of the following version of Waring's problem
(Theorem \ref{th:Waring}):

Let $F(x_1, \dots ,x_N)$ be a symmetric form of odd degree $n\ge
5$ in $N=6 \cdot 2^{n-4}$ variables, and let $q$ be any rational
number. Then the equation $F(x_i)=q$ has a rational parametric
solution depending on $2n-6$ parameters.


\section {Symmetric quintics in 6 variables} \label{sec:quin}

We begin with a simple general observation:

\begin{lem} \label{th:repr}
Every symmetric form of degree $2n+1$ in $2N\ge 6$ variables has a
zero of the form
\begin{equation}
(a_1,-a_1,\dots, a_{N-1},-a_{N-1}, 1,-1) \label{eq:prim}
\end{equation}
with $a_i\neq a_j\neq 1$.
\end{lem}

\begin{proof}
By a well-known algebraic theorem, every symmetric polynomial
$P\in k[x_1,\dots ,x_N]^{S_N}$ can be represented as a polynomial
in elementary symmetric polynomials $e_1(x_1,\dots, x_N)$, \dots ,
$e_N(x_1,\dots ,x_N)$. In particular, every symmetric form can be
represented as
$$
F=\sum_{i_1+\dots +i_k=2n+1}\alpha_{i_1,\dots ,i_k} e_{i_1}\dots
e_{i_k}
$$
where $e_{i_m}(x_1,\dots ,x_N)$ stands for the elementary
symmetric polynomial of degree $i_m$.
Each summand contains a factor of odd degree, therefore they all
vanish at any point of the form (\ref{eq:prim}).
\end{proof}

\begin{defn} \label{def:prim}
We call any point of the form (\ref{eq:prim}), as well as any its
image with respect to the action of the symmetric group, a
primitive solution of $F=0$.
\end{defn}

First consider the case of diagonal quintics.

\begin{prop} \label{prop:diag}
The equation $x_1^5+x_2^5+x_3^5+x_4^5+x_5^5+x_6^5=0$ has a
parametric solution, where $x_i$ are polynomials in two parameters
$a$ and $b$, and $x_1+\dots +x_6=0$.
\end{prop}

\begin{proof}
This equation has primitive solutions $(a,-a,b,-b,1,-1)$ and
$(-1,c,-c,d,-d,1)$. Let us look for a new solution in the form
\begin{equation}
(a\cdot t-1)^5+(-a\cdot t+c)^5+(b\cdot t-c)^5+(-b\cdot
t+d)^5+(1\cdot t-d)^5+(-1\cdot t+1)^5=0. \label{eq:5t}
\end{equation}
This gives an equation in $t$. For every $t$ we have ${x_1+\dots
+x_6=0}$. When $t$ runs over the projective line, this equation
has 5 roots: $t_1=0$, $t_2= \infty $, and 3 additional roots:
$t_3$, $t_4$, $t_5$.

The main idea is as follows: to get these 3 remaining roots $t_3$,
$t_4$, $t_5$, we have an equation in $t$ of degree $3$ obtained
from (\ref{eq:5t}):
\begin{equation}
A\cdot t^3+B\cdot t^2+C\cdot t+D=0 \label{eq:3t}
\end{equation}
where $A$, $B$, $C$, $D$ are functions in  $a$, $b$, $c$, $d$; let
us require $A=0$, $B=0$. Solving this system of two equations, we
get $c$, $d$ as functions in $a$ and $b$. Then $t=-D/C$, $D$ and
$C$ are also functions in $a$ and $b$. We have
\[
\begin{aligned}
A & = & -a^4+a^4\cdot c-b^4\cdot c+b^4\cdot d-d+1,\\
B & = & a^3-a^3\cdot c^2+b^3\cdot c^2-b^3\cdot d^2+d^2-1.
\end{aligned}
\]

If $A=0$ then $c=(d\cdot (1-b^4)+a^4-1)/(a^4-b^4)$.

If $B=0$ then $c^2=(d^2\cdot (b^3-1)+1-a^3)/(b^3-a^3)$.

We have a system of two equations, where $c$ and $d$ are unknowns.
It is easy to transform this system into one equation of degree 2
in $d$; then this equation evidently has a solution $d=1$, because
$c=d=1$ is a solution of this system, so the second root is
$d=f(a,b)$.
\end{proof}

This method can be extended to more complicated cases.

\begin{thm} \label{th:quintic}
Every symmetric quintic equation $F(x_1,\dots ,x_6)=0$ has a
parametric solution where $x_i$ are polynomials in two parameters
$a$ and $b$, and $x_1+\dots +x_6=0$.
\end{thm}

\begin{proof}
Under the assumption $x_1+\dots +x_6=0$, every quintic is
representable as
\[A_1\cdot (x_1^5+\dots+x_6^5)+A_2\cdot (x_1^3+\dots+x_6^3) \cdot (x_1^2+\dots+x_6^2).\]

Indeed, any symmetric quintic can be represented as a polynomial
in diagonal forms:
$$
\begin{aligned}
A_1(x_1^5+\dots+x_6^5) & +A_2(x_1^3+\dots+x_6^3) \cdot
(x_1^2+\dots+x_6^2)+ A_3(x_1^4+\dots+x_6^4) \cdot (x_1+\dots+x_6)
\\
& + A_4(x_1^3+\dots+x_6^3) \cdot (x_1+\dots+x_6)^2+
A_5(x_1+\dots+x_6) \cdot(x_1^2+\dots+x_6^2)^2 \\
& +A_6(x_1+\dots+x_6)^3 \cdot
(x_1^2+\dots+x_6^2)+A_7(x_1+\dots+x_6)^5,
\end{aligned}
$$
and all terms, except for the first two ones, vanish. To find a
parametric solution of the resulting quintic, we use exactly the
same method as for the diagonal quintic in Proposition
\ref{prop:diag}.

This quintic has primitive solutions
\[(a,-a,b,-b,1,-1),(-1,c,-c,d,-d,1).\]
Let us look for a new solution in the form:

$A_1\cdot ((a\cdot t-1)^5+(-a\cdot t+c)^5+(b\cdot t-c)^5+(-b\cdot
t+d)^5+(1\cdot t-d)^5+(-1\cdot t+1)^5)+A_2((a\cdot t-1)^3+(-a\cdot
t+c)^3+(b\cdot t-c)^3+(-b\cdot t+d)^3+(1\cdot t-d)^3+(-1\cdot
t+1)^3)((a\cdot t-1)^2+(-a\cdot t+c)^2+(b\cdot t-c)^2+(-b\cdot
t+d)^2+(1\cdot t-d)^2+(-1\cdot t+1)^2) = 0.$

This gives an equation in $t$. We use the same idea: this equation
in $t$ has 5 roots: $0$, $\infty$, $t_1$, $t_2$, $t_3$. For $t_1$,
$t_2$, $t_3$ we have an equation:
\[R_1\cdot t^3+R_2\cdot t^2+R_3\cdot t+R_4 =0\]
where $R_1$, $R_2$, $R_3$, $R_4$ are functions in $a$, $b$, $c$,
$d$, $A_1$, $A_2$. Now we require: $R_1=0$, $R_2=0$. This is a
system of two equations with unknowns $c$ and $d$. The equation
$R_1=0$ is linear in $c$ and $d$, and the equation $R_2=0$ is
quadratic in $c$ and $d$. It is easy to see that if we transform
this system to one equation of degree $2$ in $d$, this equation
evidently has the root $d_1=1$ corresponding to the primitive
solution
\[(a\cdot t-1,-a\cdot t+1,b\cdot t-1,-b\cdot t+1,t-1,-t+1).\]
The second root is $d_2=F(a,b)$. So the equation in $t$ reduces to
the equation: $$K_1(a,b)\cdot t+K_2(a,b)=0$$ from which we can
express $t$ as a function in $a$ and $b$. We thus find a
parametric solution for any symmetric diophantine equation of
degree 5 in 6 variables: $x_i=f_i(a,b)$ for which we have
$x_1+\dots+x_6=0.$
\end{proof}

\section {General results} \label{sec:gen}

In this section we generalize Theorem \ref{th:quintic} to the case
of a form of an arbitrary odd degree.


\begin{thm} \label{th:general}
Let $F$ be a symmetric form in $N$ variables  of odd degree $n\ge
5$ with rational coefficients. If $N\ge 6\cdot 2^{n-5}$ then the
equation $F(x_1,\dots ,x_N)=0$ has a parametric solution where
$x_i$ are polynomials in $s=2n-8$ parameters, and $x_1+\dots
+x_N=0$.
\end{thm}

\begin{proof}
First consider the equation
\begin{equation} \label{eq:III}
\begin{aligned}
(x_1^5+\dots +x_6^5) & +(x_1^4-x_2^4+x_3^4-x_4^4+x_5^4-x_6^4)\cdot
D_1(x_i) \\
& +(x_1^3+x_2^3+x_3^3+x_4^3+x_5^3+x_6^3)\cdot D_2(x_i)
\\
& +(x_1^2-x_2^2+x_3^2-x_4^2+x_5^2-x_6^2)\cdot D_3(x_i)=0,
\end{aligned}
\end{equation}
where $D_1(x_i)$ is a polynomial in $x_i$ of degree at most 1,
$D_2(x_i)$ is a polynomial in $x_i$ of degree at most 2,
$D_3(x_i)$ is a polynomial in $x_i$ of degree at most 3. It is
easy to see that this equation has solutions:
\[(a,-a,b,-b,1,-1), (-1,c,-c,d,-d,1).\]
Let us try to find a new solution by the same method:
$$x_1=a\cdot t-1, x_2=-a\cdot t+c, x_3=b\cdot t-c, x_4=-b\cdot t+d,
x_5=t-d, x_6=-t+1,$$ so equation (\ref{eq:III}) in $t$ has roots
$t=0$, $t=\infty $, and hence can be transformed to an equation of
degree $3$
\[S_1\cdot t^3+S_2\cdot t^2+S_3\cdot t+S_4=0\]
(here $S_1$, $S_2$, $S_3$, $S_4$ are functions in $a$, $b$, $c$,
$d$).

Now, as above, we have two equations for $c$ and $d$: $S_1=0$,
$S_2=0$. The equation $S_1=0$ is linear in $c$ and $d$, and
$S_2=0$ is quadratic in $c$ and $d$. From the equation $S_1=0$ it
follows that if $d=1$ then $c=1$, because for $c=1$, $d=1$,
equation (\ref{eq:III}) becomes an identity. Therefore the
equation $S_2=0$ must have the root $d_2=1$, $d_1=F(a,b)$, and
equation (\ref{eq:III}) is transformed to $S_3\cdot t+S_4=0$. So
we can set $t=-S_4/S_3$ and find a parametric solution of equation
(\ref{eq:III}).

Now suppose that we have an arbitrary symmetric form $F$ of degree
$n=2k+1$, and the number of variables is $N=4s$. For every
quadruple of variables we use a transformation of the form
\begin{equation}  \label{eq:I}
z_1=x_1+c_1, z_2=-x_2+d_1, z_3=-x_1-d_1, z_4=x_2-c_1.
\end{equation}

We represent this form as follows:

$$z_1^{2n+1}+\dots +z_N^{2n+1}+A_1(z_1^{2n-1}+\dots
+z_N^{2n-1})R_1(z_i)+\dots +A_k(z_1+\dots +z_N)R_k(z_i),$$ where
$R_j(z_i)$ are symmetric polynomials. Let us now look how this
transformation works.

It is easy to see that we obtain a form of degree $2n$ whose
coefficients are functions in $c_N$, $d_N$, but we obtain a new
symmetric construction of the form
\begin{equation} \label{eq:II}
x_1^{2n}-x_2^{2n}+x_3^{2n}-x_4^{2n}+\dots
\end{equation}

Under transformation (\ref{eq:I}), polynomials of this type go
over to linear combinations of symmetric diagonal polynomials of
odd degree and polynomials of the form  (\ref{eq:II}) of even
degree. The degree of these polynomials is smaller than $2n+1$.
Repeating this transformation several times, we obtain an equation
of the form
$$
\begin{aligned}
(x_1^5+\dots +x_6^5) & +(x_1^4-x_2^4+x_3^4-x_4^4+x_5^4-x_6^4)\cdot
D_1(x_i) \\
& +(x_1^3+x_2^3+x_3^3+x_4^3+x_5^3+x_6^3)\cdot D_2(x_i)
\\
& +(x_1^2-x_2^2+x_3^2-x_4^2+x_5^2-x_6^2)\cdot D_3(x_i)=0,
\end{aligned}
$$
but for this equation the existence of a solution depending on two
parameters $a$ and $b$ is proved in Theorem \ref{th:quintic}. So
if a symmetric form has $6\cdot 2^{n-5}$ ($n\ge 5$ is an odd
number) variables, it has a parametric solution in $s=2n-8$
parameters.
\end{proof}

\medskip

One can use the same approach as in Theorem \ref{th:general} for
finding parametric solutions of the following problem of Waring's
type.

\begin{thm} \label{th:Waring}
Let $n\ge 5$ be an odd integer, $N=6\cdot 2^{n-4}$, let $F$ be a
symmetric form of degree $n$ with rational coefficients, and let
$q$ be a fixed rational number. Then the equation $F(x_1,\dots
,x_N)=q$ has a parametric solution depending on $s=2n-6$
parameters.
\end{thm}

\begin{proof} The proof follows, almost word for word, the
arguments of Theorem \ref{th:general}, and we leave the details to
a scrupulous reader.
\end{proof}

\begin{remark} \label{th:repr}
One can pose the following problem: for given integers $n$ and $k$
to find an integer $N$ such that every symmetric form in $N$
variables of degree $2n+1$ with rational coefficients has a
parametric solution depending on $k$ parameters. Such a problem
can be treated using the approach proposed in this paper. In
particular, it is very easy to generalize the results of this
paper from symmetric quintics  in 6 variables to symmetric
quintics in any even number of variables. One can use primitive
solutions to construct, for every symmetric quintic in $2P$
variables, a parametric solution depending on $P-1$ parameters.
Using transformation (\ref{eq:I}), one can then generalize this
result for every symmetric form of odd degree. Finally, for every
given integers $n$ and $k$ one can find an integer $N$ such that
every symmetric form in $N$ variables of degree $2n+1$ with
rational coefficients has a parametric solution depending on $s$
parameters with $s>k$. If in this solution we substitute arbitrary
rational numbers instead of $s-k$ parameters, we obtain a desired
parametric solution depending on $k$ parameters.
\end{remark}

\end{document}